\documentclass[a4paper, 12pt,oneside,reqno]{amsart}
\usepackage{graphicx} 
\usepackage{tikz}

\begin{document}

\title[Gauss diagrams as cubic graphs]{Gauss diagrams as cubic graphs: The choice of the Hamiltonian cycle matters}

\author{Alexei Vernitski}

\address{University of Essex, Colchester, UK}

\maketitle

\begin{abstract}
We explore to what extent the properties of a Gauss diagram are affected by the choice of its Hamiltonian cycle. We present an example of a realizable Gauss diagram and an unrealizable Gauss diagram that differ only by a choice of the Hamiltonian cycle.  We present an example of two Gauss diagrams that correspond to different curves and differ only by a choice of the Hamiltonian cycle. We prove that a certain natural type of change of the Hamiltonian cycle preserves the realizability of the Gauss diagram.
\end{abstract}



\section{Introduction}

Gauss diagrams are an important instrument in knot theory. In graph-theoretical terms, a Gauss diagram is a pair consisting of a cubic bipartite Hamiltionian graph and a Hamiltonian cycle in this graph. Traditionally, when one draws a Gauss diagram, the fixed Hamiltonian cycle forms the circumference of a circle, and the other edges are drawn as chords inside the circle. Please see three examples in Figures \ref{fig:Gauss-diagram-1}, \ref{fig:Gauss-diagram-2}, \ref{fig:Gauss-diagram-3}; these diagrams not only illustrate the traditional shape of Gauss diagrams, but pertain to examples given below.

\begin{figure}
\centering 
\begin{minipage}{.33\textwidth}
\centering 
\begin{tikzpicture} [thick, scale=0.5]
\draw[fill=none, line width = 2](0,0) circle (3);

\draw (0*36:3) -- (3*36:3);
\draw (2*36:3) -- (5*36:3);
\draw (4*36:3) -- (7*36:3);
\draw (6*36:3) -- (9*36:3);
\draw (8*36:3) -- (11*36:3);

\end{tikzpicture}
\caption{}
\label{fig:Gauss-diagram-1}
\end{minipage}%
\begin{minipage}{.33\textwidth}
\centering 
\begin{tikzpicture} [thick, scale=0.5]
\draw[fill=none, line width = 2](0,0) circle (3);

\draw (0*36:3) -- (5*36:3);
\draw (2*36:3) -- (7*36:3);
\draw (4*36:3) -- (9*36:3);
\draw (6*36:3) -- (11*36:3);
\draw (8*36:3) -- (13*36:3);

\end{tikzpicture}
\caption{}
\label{fig:Gauss-diagram-2}
\end{minipage}
\begin{minipage}{.33\textwidth}
\centering 
\begin{tikzpicture} [thick, scale=0.5]
\draw[fill=none, line width = 2](0,0) circle (3);

\draw (0*36:3) -- (5*36:3);
\draw (6*36:3) -- (9*36:3);
\draw (4*36:3) -- (1*36:3);
\draw (2*36:3) -- (7*36:3);
\draw (8*36:3) -- (3*36:3);

\end{tikzpicture}
\caption{}
\label{fig:Gauss-diagram-3}
\end{minipage}
\end{figure}

The aim of this short note is to explore to what extent the properties of a Gauss diagram are affected by the choice of its Hamiltonian cycle. The example in Section \ref{sec:realizable} demonstrates that a realizable Gauss diagram and an unrealizable Gauss diagram can share the same graph and differ only by a choice of the Hamiltonian cycle.  The example in Section \ref{sec:curve} demonstrates that two realizable Gauss diagrams that correspond to different curves can share the same graph and differ only by a choice of the Hamiltonian cycle. These examples are minimal in size. Then in Section \ref{sec:preserved} we prove that a certain natural type of change of the Hamiltonian cycle preserves the realizability of the Gauss diagram.

\section{The choice of the Hamiltonian cycle affects realizability} \label{sec:realizable}

The graph shown in Figure \ref{fig:Mobius-ladder} is known as the M\"{o}bius ladder of order 5 or the 5-M\"{o}bius ladder. 

\begin{figure}
\centering 
\begin{tikzpicture}
    \def\squareSize{1.5}

    \draw (0, 0) -- (\squareSize, 0);
    \draw (\squareSize, 0) -- (\squareSize, \squareSize);
    \draw (\squareSize, \squareSize) -- (0, \squareSize);
    \draw (0, \squareSize) -- (0, 0);
    
    \draw (\squareSize, 0) -- (2 * \squareSize, 0);
    \draw (2 * \squareSize, 0) -- (2 * \squareSize, \squareSize);
    \draw (2 * \squareSize, \squareSize) -- (\squareSize, \squareSize);
    
    \draw (2 * \squareSize, 0) -- (3 * \squareSize, 0);
    \draw (3 * \squareSize, 0) -- (3 * \squareSize, \squareSize);
    \draw (3 * \squareSize, \squareSize) -- (2 * \squareSize, \squareSize);
    
    \draw (3 * \squareSize, 0) -- (4 * \squareSize, 0);
    \draw (4 * \squareSize, 0) -- (4 * \squareSize, \squareSize);
    \draw (4 * \squareSize, \squareSize) -- (3 * \squareSize, \squareSize);
    
    \draw plot [smooth,tension=0.8] coordinates {(0, \squareSize) (\squareSize, 1.5 *\squareSize) (3*\squareSize, 1.5 *\squareSize) (3.5 * \squareSize + \squareSize, \squareSize) (3 * \squareSize + \squareSize, 0)};
    
    \draw plot [smooth,tension=0.8] coordinates {(0, 0) (\squareSize, -0.5 *\squareSize) (3*\squareSize, -0.5 *\squareSize) (3.5 * \squareSize + \squareSize, 0) (3 * \squareSize + \squareSize, \squareSize)};    

\end{tikzpicture}
\caption{The M\"{o}bius ladder of order 5}
\label{fig:Mobius-ladder}
\end{figure}

If we choose a Hamiltonian cycle in the 5-M\"{o}bius ladder as highlighted in Figure \ref{fig:Mobius-ladder-1}, this combination of a graph and a Hamiltonian cycle in the graph is equivalent to the Gauss diagram shown in Figure \ref{fig:Gauss-diagram-1}. It is known (and easy to check) that this Gauss diagram is not realizable, that is, it does not correspond to a closed curve on the plane.

\begin{figure}
\centering 
\begin{tikzpicture}
    \def\squareSize{1.5}

    \draw[line width = 2] (0, 0) -- (\squareSize, 0);
    \draw[line width = 2] (\squareSize, 0) -- (\squareSize, \squareSize);
    \draw (\squareSize, \squareSize) -- (0, \squareSize);
    \draw[line width = 2] (0, \squareSize) -- (0, 0);
    
    \draw (\squareSize, 0) -- (2 * \squareSize, 0);
    \draw[line width = 2] (2 * \squareSize, 0) -- (2 * \squareSize, \squareSize);
    \draw[line width = 2] (2 * \squareSize, \squareSize) -- (\squareSize, \squareSize);
    
    \draw[line width = 2] (2 * \squareSize, 0) -- (3 * \squareSize, 0);
    \draw[line width = 2] (3 * \squareSize, 0) -- (3 * \squareSize, \squareSize);
    \draw (3 * \squareSize, \squareSize) -- (2 * \squareSize, \squareSize);
    
    \draw (3 * \squareSize, 0) -- (4 * \squareSize, 0);
    \draw[line width = 2] (4 * \squareSize, 0) -- (4 * \squareSize, \squareSize);
    \draw[line width = 2] (4 * \squareSize, \squareSize) -- (3 * \squareSize, \squareSize);
    
    \draw[line width = 2] plot [smooth,tension=0.8] coordinates {(0, \squareSize) (\squareSize, 1.5 *\squareSize) (3*\squareSize, 1.5 *\squareSize) (3.5 * \squareSize + \squareSize, \squareSize) (3 * \squareSize + \squareSize, 0)};
    
    \draw plot [smooth,tension=0.8] coordinates {(0, 0) (\squareSize, -0.5 *\squareSize) (3*\squareSize, -0.5 *\squareSize) (3.5 * \squareSize + \squareSize, 0) (3 * \squareSize + \squareSize, \squareSize)};    

\end{tikzpicture}
\caption{Producing the Gauss diagram from Figure \ref{fig:Gauss-diagram-1}}
\label{fig:Mobius-ladder-1}
\end{figure}

If we choose a Hamiltonian cycle in the 5-M\"{o}bius ladder as highlighted in Figure \ref{fig:Mobius-ladder-2}, this combination of a graph and a Hamiltonian cycle in the graph is equivalent to the Gauss diagram shown in Figure \ref{fig:Gauss-diagram-2}. This Gauss diagram is realizable; indeed, it corresponds to the standard diagram of the torus knot $5_1$, see Figure \ref{fig:5-1}.

\begin{figure}
\centering 
\begin{tikzpicture}
    \def\squareSize{1.5}

    \draw[line width = 2] (0, 0) -- (\squareSize, 0);
    \draw (\squareSize, 0) -- (\squareSize, \squareSize);
    \draw[line width = 2] (\squareSize, \squareSize) -- (0, \squareSize);
    \draw (0, \squareSize) -- (0, 0);
    
    \draw[line width = 2] (\squareSize, 0) -- (2 * \squareSize, 0);
    \draw (2 * \squareSize, 0) -- (2 * \squareSize, \squareSize);
    \draw[line width = 2] (2 * \squareSize, \squareSize) -- (\squareSize, \squareSize);
    
    \draw[line width = 2] (2 * \squareSize, 0) -- (3 * \squareSize, 0);
    \draw (3 * \squareSize, 0) -- (3 * \squareSize, \squareSize);
    \draw[line width = 2] (3 * \squareSize, \squareSize) -- (2 * \squareSize, \squareSize);
    
    \draw[line width = 2] (3 * \squareSize, 0) -- (4 * \squareSize, 0);
    \draw (4 * \squareSize, 0) -- (4 * \squareSize, \squareSize);
    \draw[line width = 2] (4 * \squareSize, \squareSize) -- (3 * \squareSize, \squareSize);
    
    \draw[line width = 2] plot [smooth,tension=0.8] coordinates {(0, \squareSize) (\squareSize, 1.5 *\squareSize) (3*\squareSize, 1.5 *\squareSize) (3.5 * \squareSize + \squareSize, \squareSize) (3 * \squareSize + \squareSize, 0)};
    
    \draw[line width = 2] plot [smooth,tension=0.8] coordinates {(0, 0) (\squareSize, -0.5 *\squareSize) (3*\squareSize, -0.5 *\squareSize) (3.5 * \squareSize + \squareSize, 0) (3 * \squareSize + \squareSize, \squareSize)};    

\end{tikzpicture}
\caption{Producing the Gauss diagram from Figure \ref{fig:Gauss-diagram-2}}
\label{fig:Mobius-ladder-2}
\end{figure}

\begin{figure}
\centering 
\begin{minipage}{.5\textwidth}
\centering 
\begin{tikzpicture}
\draw[thick] (-2,1) to[out=120, in=150] (0,3);
\draw[thick] (-2,1) to[out=30, in=210] (0,3);

\draw[thick] (2,1) to[out=60, in=30] (0,3);
\draw[thick] (2,1) to[out=160, in=330] (0,3);

\draw[thick] (-2,1) to[out=180+120, in=150] (-1,0);
\draw[thick] (-2,1) to[out=180+30, in=210] (-1,0);

\draw[thick] (2,1) to[out=60+180, in=30] (1,0);
\draw[thick] (2,1) to[out=160+180, in=330] (1,0);

\draw[thick] (-1,0) to[out=180+150, in=180+30] (1,0);
\draw[thick] (-1,0) to[out=180+210, in=330+180] (1,0);

\end{tikzpicture}
\caption{Curve $5_1$}
\label{fig:5-1}
\end{minipage}%
\begin{minipage}{.5\textwidth}
\centering 
\begin{tikzpicture}
\draw[thick] (-2,1) to[out=120, in=150] (0,3);
\draw[thick] (-2,1) to[out=30, in=210] (0,3);

\draw[thick] (2,1) to[out=60, in=30] (0,3);
\draw[thick] (2,1) to[out=160, in=330] (0,3);

\draw[thick] (-2,1) to[out=180+120, in=150] (0,1);
\draw[thick] (-2,1) to[out=180+30, in=210] (0,0);

\draw[thick] (2,1) to[out=60+180, in=30] (0,1);
\draw[thick] (2,1) to[out=160+180, in=330] (0,0);

\draw[thick] (0,1) to[out=180+150, in=180+210] (0,0);
\draw[thick] (0,1) to[out=180+30, in=330+180] (0,0);

\end{tikzpicture}
\caption{Curve $5_2$}
\label{fig:5-2}
\end{minipage}
\end{figure}

Comparing examples presented in Figures \ref{fig:Mobius-ladder-1} and \ref{fig:Mobius-ladder-2}, we conclude that a realizable Gauss diagram and an unrealizable Gauss diagram can share the same cubic graph and differ only by a choice of the Hamiltonian cycle.

\section{The choice of the Hamiltonian cycle affects the curve} \label{sec:curve}

\begin{figure}
\centering 
\begin{tikzpicture}
    \def\squareSize{1.5}

    \draw[line width = 2] (0, 0) -- (\squareSize, 0);
    \draw (\squareSize, 0) -- (\squareSize, \squareSize);
    \draw[line width = 2] (\squareSize, \squareSize) -- (0, \squareSize);
    \draw[line width = 2] (0, \squareSize) -- (0, 0);
    
    \draw[line width = 2] (\squareSize, 0) -- (2 * \squareSize, 0);
    \draw (2 * \squareSize, 0) -- (2 * \squareSize, \squareSize);
    \draw[line width = 2] (2 * \squareSize, \squareSize) -- (\squareSize, \squareSize);
    
    \draw[line width = 2] (2 * \squareSize, 0) -- (3 * \squareSize, 0);
    \draw (3 * \squareSize, 0) -- (3 * \squareSize, \squareSize);
    \draw[line width = 2] (3 * \squareSize, \squareSize) -- (2 * \squareSize, \squareSize);
    
    \draw[line width = 2] (3 * \squareSize, 0) -- (4 * \squareSize, 0);
    \draw[line width = 2] (4 * \squareSize, 0) -- (4 * \squareSize, \squareSize);
    \draw[line width = 2] (4 * \squareSize, \squareSize) -- (3 * \squareSize, \squareSize);
    
    \draw plot [smooth,tension=0.8] coordinates {(0, \squareSize) (\squareSize, 1.5 *\squareSize) (3*\squareSize, 1.5 *\squareSize) (3.5 * \squareSize + \squareSize, \squareSize) (3 * \squareSize + \squareSize, 0)};
    
    \draw plot [smooth,tension=0.8] coordinates {(0, 0) (\squareSize, -0.5 *\squareSize) (3*\squareSize, -0.5 *\squareSize) (3.5 * \squareSize + \squareSize, 0) (3 * \squareSize + \squareSize, \squareSize)};    

\end{tikzpicture}
\caption{Producing the Gauss diagram from Figure \ref{fig:Gauss-diagram-3}}
\label{fig:Mobius-ladder-3}
\end{figure}
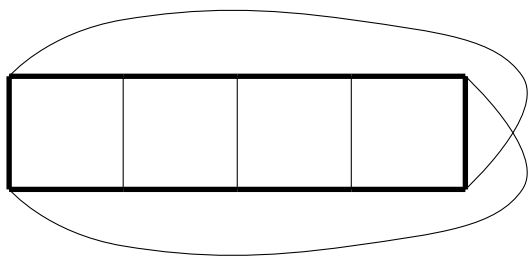

If we choose a Hamiltonian cycle in the 5-M\"{o}bius ladder as highlighted in Figure \ref{fig:Mobius-ladder-3}, this combination of a graph and a Hamiltonian cycle in the graph is equivalent to the Gauss diagram shown in Figure \ref{fig:Gauss-diagram-3}. This Gauss diagram is realizable; indeed, it corresponds to the standard diagram of the knot $5_2$, see Figure \ref{fig:5-2}.

Comparing examples presented in Figures \ref{fig:Mobius-ladder-2} and \ref{fig:Mobius-ladder-3}, we conclude that two realizable Gauss diagrams that correspond to different curves (in this example, the curves in Figures \ref{fig:5-1} and \ref{fig:5-2}) can share the same graph and differ only by a choice of the Hamiltonian cycle.

\section{When realizability is preserved} \label{sec:preserved}

After seeing examples in Sections \ref{sec:realizable} and \ref{sec:curve} showing that the change of the Hamiltonian cycle alters the properties of a Gauss diagram drastically, you might wonder if there are any contexts in which the change of the Hamiltonian path preserves some properties of the Gauss diagram. As we will demonstrate now, there is a natural type of change of the Hamiltonian cycle that preserves the realizability of the Gauss diagram.

Consider a Hamiltonian cycle containing the vertices of the graph in the order $v_1, \dots, v_n$, and reverse the ordering of vertices in a fragment of this cycle (say, $v_{k+1}, \dots, v_n$), producing a new ordering of vertices $v_1, \dots, v_k, v_n, v_{n-1}, \dots, v_{k+1}$. If the latter sequence of vertices $v_1, \dots, v_k, v_n, v_{n-1}, \dots, v_{k+1}$ is also a Hamiltonian cycle, we shall say that these two Hamiltonian cycles are produced from each other by a \textit{single flip}. In the language of Gauss diagrams, consider a Gauss diagram shown in Figure \ref{fig:flip-1} consisting of two fragments $F1$ and $F2$ and two intersecting chords in adjacent positions. A single flip applied between $F1$ and $F2$ changes the Hamilton path as highlighted in Figure \ref{fig:flip-2}. Another way of describing the single flip applied between $F1$ and $F2$ is to flip $F2$ (together with all chord ends attached to $F2$) as shown in Figure \ref{fig:flip-3}. If you inspect the two Gauss diagrams in Figures \ref{fig:Gauss-diagram-2} and \ref{fig:Gauss-diagram-3}, you will see that they are produced from one another by a single flip.

\begin{figure}
\centering 
\begin{minipage}{.33\textwidth}
\centering 
\begin{tikzpicture} [thick, scale=0.5]
\draw[fill=none, line width = 2](0,0) circle (3);

\draw (-2*36:3) -- (3*36:3);
\draw (2*36:3) -- (-3*36:3);

\filldraw[fill=white] (-3.5, -1.5) rectangle (-1.5, 1.5);
\node at (-2.5, 0) {$F1$};

\filldraw[fill=white] (1.5, -1.5) rectangle (3.5, 1.5);
\node at (2.5, 0) {$F2$};

\end{tikzpicture}
\caption{}
\label{fig:flip-1}
\end{minipage}%
\begin{minipage}{.33\textwidth}
\centering 
\begin{tikzpicture} [thick, scale=0.5]
\draw[fill=none](0,0) circle (3);
\draw[line width = 2] (0,0) ++(36*3:3) arc (36*3:36*7:3);
\draw[line width = 2] (0,0) ++(36*2:3) arc (36*2:0:3);
\draw[line width = 2] (0,0) ++(0:3) arc (0:-36*2:3);

\draw[line width = 2] (-2*36:3) -- (3*36:3);
\draw[line width = 2] (2*36:3) -- (-3*36:3);

\filldraw[fill=white] (-3.5, -1.5) rectangle (-1.5, 1.5);
\node at (-2.5, 0) {$F1$};

\filldraw[fill=white] (1.5, -1.5) rectangle (3.5, 1.5);
\node at (2.5, 0) {$F2$};

\end{tikzpicture}
\caption{}
\label{fig:flip-2}
\end{minipage}
\begin{minipage}{.33\textwidth}
\centering 
\begin{tikzpicture} [thick, scale=0.5]
\draw[fill=none, line width = 2](0,0) circle (3);

\draw (-2*36:3) -- (3*36:3);
\draw (2*36:3) -- (-3*36:3);

\filldraw[fill=white] (-3.5, -1.5) rectangle (-1.5, 1.5);
\node at (-2.5, 0) {$F1$};

\filldraw[fill=white] (1.5, -1.5) rectangle (3.5, 1.5);
\node[yscale=-1] at (2.5, 0) {$F2$};

\end{tikzpicture}
\caption{}
\label{fig:flip-3}
\end{minipage}
\end{figure}

Theorem. Changing the Hamiltonian cycle by applying a single flip preserves the realizability of a Gauss diagram. 

\begin{proof}
Consider a realizable Gauss diagram $(G, H)$ consisting of a graph $G$ and a Hamiltonian cycle $H$. Since $(G, H)$ is realizable, there is a closed curve $C$ corresponding to $(G, H)$. Suppose a Hamiltonian cycle $H'$ is produced from $H$ by a single flip. We will show that we can amend $C$ in such a way that the Gauss diagram of the changed curve $C'$ is $(G, H')$; hence, $(G, H')$ is realizable. 

Indeed, suppose the single flip which turns $H$ into $H'$ is applied to a two-chord fragment of $(G, H)$ which looks as in Figure \ref{fig:proof-1}; what is not shown in this diagram is a part of $H$ connecting $a$ with $c$ and a part of $H$ connecting $b$ with $d$. In curve $C$, this fragment of the Gauss diagram corresponds to a fragment of the curve shown in Figure \ref{fig:proof-2}; what is not shown in this diagram is a part of $C$ connecting $a$ with $c$ and a part of $C$ connecting $b$ with $d$. Produce curve $C'$ by replacing the fragment shown in Figure \ref{fig:proof-2} by the fragment shown in Figure \ref{fig:proof-3}. (An example of this kind of change can be seen in the two curves shown in Figures \ref{fig:5-1} and \ref{fig:5-2}, with two crossings altered in the middle of the bottom part of the curves.) Outside this fragment there are a part of $C'$ connecting $a$ with $c$ and a part of $C'$ connecting $b$ with $d$, therefore, $C'$ is one closed curve. It is easy to see that the Gauss diagram of $C'$ is $(G, H')$.
\end{proof}

\begin{figure}
\centering 
\begin{minipage}{.33\textwidth}
\centering 
\begin{tikzpicture} [thick, scale=0.5]
\draw[line width = 2] (0,0) ++(36*3+10:3) arc (36*3+10:36*2-10:3);
\draw[line width = 2] (0,0) ++(-36*2+10:3) arc (-36*2+10:-36*3-10:3);

\draw (-2*36:3) -- (3*36:3);
\draw (2*36:3) -- (-3*36:3);

\node at (-2.2, 3) {$a$};
\node at (2.2, 3) {$b$};
\node at (-2.2, -3) {$c$};
\node at (2.2, -3) {$d$};

\end{tikzpicture}
\caption{}
\label{fig:proof-1}
\end{minipage}%
\begin{minipage}{.33\textwidth}
\centering 
\begin{tikzpicture} [thick, scale=0.5]

\draw plot [smooth,tension=0.8] coordinates {(0, 0.5) (3, 5) (6, 0.5)};
\draw plot [smooth,tension=0.8] coordinates {(0, 6-0.5) (3, 6-5) (6, 6-0.5)};

\node at (0, 6) {$a$};
\node at (6, 6) {$b$};
\node at (0, 0) {$d$};
\node at (6, 0) {$c$};

\end{tikzpicture}
\caption{}
\label{fig:proof-2}
\end{minipage}
\begin{minipage}{.33\textwidth}
\centering 
\begin{tikzpicture} [thick, scale=0.5]
\draw plot [smooth,tension=0.8] coordinates {(0.5, 0) (5, 3) (0.5, 6)};
\draw plot [smooth,tension=0.8] coordinates {(6-0.5, 0) (6-5, 3) (6-0.5, 6)};

\node at (0, 6) {$a$};
\node at (6, 6) {$b$};
\node at (0, 0) {$d$};
\node at (6, 0) {$c$};

\end{tikzpicture}
\caption{}
\label{fig:proof-3}
\end{minipage}
\end{figure}

\section{Context}

Recall that a Gauss diagram is a pair $(G, H)$ consisting of a cubic bipartite Hamiltonian graph $G$ and a Hamiltonian cycle $H$ in $G$. One could think of $G$ as a summary of $(G, H)$, and as we saw in Sections \ref{sec:realizable}, \ref{sec:curve}, this summary does not retain much interesting information about the Gauss diagram. Another graph-theoretical construction arising from Gauss diagrams is so-called interlacement graph of the Gauss diagram. In this graph, the vertices correspond to the chords of the Gauss diagram, and two vertices are connected with an edge if and only if the two chords intersect in the Gauss diagram. It is known that a Gauss diagram is realizable if and only if its interlacement graph satisfies certain conditions \cite{rosenstiehl:hal-00259712,10.1007/3-540-63938-1_65,DBLP:journals/dm/ShtyllaTZ09,doi:10.1142/S0218216523500591}.

Several recent papers considered moves applied to Gauss diagrams; usually this is done in the context of representing Reidemeister moves in the language of Gauss diagrams \cite{nelson2001unknotting,manturov2010parity,mortier2013polyak,boden2015bridge,Valette16,kaur2017gauss,meilhan2019arrow,ito2021gauss,xue2023unknotting}. It seems that the move we study in Section \ref{sec:preserved}, which we call the single flip, has not been studied. Having said this, there is a similarly-looking move which forms a part of some algorithms for detecting if a Gauss diagram is realizable \cite{rosenstiehl1984gauss,grinblat2020realizabilty} (and which also features in other contexts under a name of Conway's smoothing of crossings); that move preserves the realizability of Gauss diagrams.

\bibliographystyle{elsarticle-num}
\bibliography{bib}

\end{document}